\documentclass[12pt]{amsart}
\usepackage{epsfig,amsfonts,amsthm,amssymb,latexsym,amsmath}
\usepackage[latin1]{inputenc}
\usepackage{bm}
\usepackage{cite}

\newfont{\bms}{msbm10 scaled 1100}
\def\matR{\mbox{\bms R}}
\def\matQ{\mbox{\bms Q}}
\def\matZ{\mbox{\bms Z}}
\def\matK{\mbox{\bms K}}
\def\matL{\mbox{\bms L}}
\def\matC{\mbox{\bms C}}

\newfont{\bmsi}{msbm10 scaled 700}

\def\mattQ{\mbox{\bmsi Q}}
\def\mattZ{\mbox{\bmsi Z}}
\def\mattK{\mbox{\bmsi K}}

\newfont{\bmsw}{msbm10 scaled 2400}

\newfont{\bmst}{msbm10 scaled 1600}

\newcommand{\bb}{\hfill \rule{2mm}{2mm}} 

\theoremstyle{plain}
\newtheorem{theorem}{Theorem}[section]
\newtheorem{corollary}[theorem]{Corollary}
\newtheorem{remark}[theorem]{Remark}
\theoremstyle{definition}
\newtheorem{definition}[theorem]{Definition}
\theoremstyle{example}

\theoremstyle{proposition}
\newtheorem{proposition}[theorem]{Proposition}
\theoremstyle{lemma}

\vspace{5pt}

\title{Rotated $D_n$-lattices}

\vspace{5pt}

\thanks{This work was partially supported by  $^{1}$(CNPq 140239/2009-0), $^{2}$(CNPq   569966/2008-6), $^{3}$(CNPq 309561/2009-4) and FAPESP 2007/56052-8}

\author{\scriptsize Grasiele C. \ Jorge$^{1}$, \  Agnaldo J. \ Ferrari$^{2}$ \ and \
\ Sueli I. R. \ Costa$^{3}$}
\date{}

\begin{document}
\maketitle

\vspace{-20pt}
\begin{center}
{\footnotesize $^{1,3}$ UNICAMP - Universidade Estadual de Campinas,
13083-859, Campinas, SP, BRAZIL \\ $^{2}$ UFLA - Universidade
Federal de Lavras,
37200-000, Lavras, MG, BRAZIL \\
Email addresses: [grajorge, ferrari, sueli] \ @ime.unicamp.br \\
}\end{center}

\hrule

\begin{abstract}
Based on algebraic number theory we construct some families of
rotated $D_n$-lattices with full diversity which can be good for
signal transmission over both Gaussian and Rayleigh \linebreak
fading channels. Closed-form expressions for the minimum  product
distance of those lattices are obtained through  algebraic
properties.
\end{abstract}

\medskip
\noindent 
\medskip

\hrule

\section{Introduction}

A lattice $\Lambda=\Lambda^{n} \subseteq \mathbb{R}^n$ is a discrete
set generated by integer combinations of $n$ linearly independents
vectors ${\bm v_1},\ldots,{\bm v_n} \in \mathbb{R}^n$.  Its packing
density $\Delta(\Lambda)$  is the proportion of the space
$\matR^{n}$ covered by congruent disjoint spheres of maximum radius
\cite{sloane}. A lattice $\Lambda$ has diversity $m\leq n$ if $m$ is
the maximum number such that for all  ${\bm y}=(y_1,\cdots,y_n) \in
\Lambda$, ${\bm y} \neq {\bm 0}$ there are at least $m$
non-vanishing coordinates. Given a full diversity lattice $\Lambda
\subseteq \matR^{n}$ $(m=n)$, the minimum product distance  is
defined as $d_{min}(\Lambda) = \min\{\prod_{i=1}^{n}|y_i|\,\,
\mbox{for all}\,\, {\bm y}=(y_1,\cdots,y_n) \in \Lambda,{\bm  y}
\neq {\bm 0} \}$ \cite{Oggier}.

Signal constellations having lattice structure have been studied as
meaningful means for \linebreak signal transmission over both
Gaussian and single-antenna Rayleigh  fading channel \cite{boutros}.
\linebreak Usually the problem of finding good signal constellations
for a Gaussian channel is  associated to the search for lattices
with high packing density \cite{sloane}. On the other hand for a
Rayleigh fading channel the efficiency,  measured by lower error
probability in the  transmission, is strongly related to the lattice
diversity and minimum  product distance \cite{boutros},
\cite{Oggier}. The approach in this work, following \cite{eva2} and
\cite{Oggier} is the use of  algebraic number theory to construct
lattices with good  performance for both channels.

For general lattices the packing density and the minimum product
distance are usually hard to estimate \cite{mic}.  Those parameters
can be obtained in certain cases of   lattices associated to number
fields, through algebraic properties.

In \cite{gabriele}, \cite{Oggier} and \cite{upper-bound} some
families of rotated $\matZ^{n}$-lattices with full diversity and
good minimum product distance are studied for transmission over
Rayleigh fading channels. In \cite{suarez} the lattices $A_{p-1}$,
$p$ prime, $E_6$, $E_8$ $K_{12}$ and $\Lambda_{24}$ were realized as
full diversity ideal lattices via some subfields of cyclotomic
fields.  In \cite{boutros} rotated $n$-dimensional lattices
(including $D_4$, $K_{12}$ and $\Lambda_{16}$) which are good for
both channels are constructed with diversity $n/2$.

In this work we also attempt to consider lattices which are feasible
for both channels by constructing rotated $D_n$-lattices with full
diversity $n$ and get a closed-form for their minimum product
distance. The results were obtained for $n=2^{r-2},$ $r \geq 5$ and
$n=(p-1)/2$, $p$ prime and $p\geq 7,$ in Propositions \ref{Prop1},
\ref{idealI} and \ref{rotacionado}. As it is known, a $D_n$ lattice
has better packing density $\delta(D_n)$  when compared to
$\matZ^{n}$ ($D_n$ has the best lattice packing density for
$n=3,4,5$ and $\lim_{n\longrightarrow
\infty}\frac{\delta(\mattZ^{n})}{\delta(D_n)} =0$) and also a very
efficient decoding algorithm \cite{sloane}. The relative minimum
product distances $d_{p,rel}(D_n)$ of the rotated $D_n$-lattices
obtained here are  smaller than the  minimum product distance
$d_{p,rel}(\matZ^{n})$ of rotated $\matZ^{n}$-lattices constructed
for the Rayleigh channels in \cite{and} and \cite{Oggier}, but, as
it is shown  in Sections 4 and 5, $\lim_{n\longrightarrow
\infty}\frac{\sqrt[n]{d_{p,rel}(\mattZ^{n})}}{\sqrt[n]{d_{p,rel}(D_n)}}
=\sqrt{2}$, what offers a good trade-off.

In Sections 2 e 3 we summarize some definitions and results on
Algebraic Number Theory. Sections 4 and 5 are devoted to the
construction of full diversity rotated $D_n$-lattices through
cyclotomic fields and the deduction of their minimum product
distance.

\section{Number Fields}

In this section we summarize some concepts and results of algebraic
\linebreak number theory and establish the notation to be used from
now on. The results presented here can be found in  \cite{traj},
\cite{pierre}, \cite{stewart} and \cite{was}.

Let $\matK$ be a number field of degree $n$ and $\mathcal
O_{\mattK}$ its ring of integers. It can be shown that every nonzero
fractionary ideal $I$ of $\mathcal O_{\mattK}$ is a free
$\matZ$-module of rank $n$.

There are exactly $n$ distinct $\matQ$-homomorphisms
$\{\sigma_i\}_{i=1}^{n}$ of $\matK$ in $\matC.$ A homomorphism
$\sigma_i$ is said {\it real} if $\sigma_i(\matK)\subset \matR$, and
the field $\matK$ is said {\it totally real} if $\sigma_i$ is real
for all $i=1,\cdots,n.$

Given $x \in \matK,$ the values $N(x) = N_{\mattK|\mattQ}(x) =
\prod_{i=1}^{n} \sigma_i(x), \ \ \ Tr(x)=Tr_{\mattK|\mattQ}(x) =
\sum_{i=1}^{n} \sigma_i(x)$ are called, {\it norm} and {\it trace}
of $x$ in $\matK|\matQ,$  respectively. It can shown that if $x \in
\mathcal{O}_{\mattK}$, then $N(x), Tr(x) \in \matZ$.

Let $\{\omega_1,\ldots,\omega_n\}$ be a $\matZ$-basis of $\mathcal
O_{\mattK}.$ The integer  $d_{\mattK} =  (det[
\sigma_j(\omega_i)]_{i,j=1}^{n})^2$ is called the {\it discriminant}
of $\matK$.

The {\it norm} of an ideal $I \subseteq \mathcal O_{\mattK}$  is
defined as $N(I) = |\mathcal O_{\mattK}/I|.$

The {\it codifferent} de $\matK|\matQ$ is the fractionary ideal
$\Delta(\matK|\matQ)^{-1}=\{x\in \matK; \,\, \forall \, \alpha \in
{\mathcal{O}}_{\mattK},\,Tr_{\mattK|\mattQ}(x\alpha) \in \matZ\}$ of
${\mathcal{O}}_{\mattK}.$

Let $\zeta=\zeta_m \in \matC$ be a  primitive $m$-th root of unity.
We consider here the {\it cyclotomic field} $\matQ(\zeta)$ and its
subfield $\matK=\matQ(\zeta+\zeta^{-1}).$ We have that
$[\matQ(\zeta+\zeta^{-1}):\matQ]= \varphi(m)/2$, where $\varphi$ is
the Euler function;  $\mathcal O_{\mattK} = \matZ
[\zeta+\zeta^{-1}];$ $d_{\mattK}=p^{\frac{p-3}{2}}$  if $m=p$, $p$
prime, $p\geq 5$ and $d_{\mattK}=2^{(r-1)2^{r-2}-1}$ if $m=2^r$.

\section{Ideal lattices}

 From now on, let $\matK$ be a totally real number
field. Let $\alpha \in \matK$ such that $\alpha_i=\sigma_i(\alpha)
> 0$ for all $i=1,\cdots,n.$ The homomorphism

$$\left.\begin{array}{l} \sigma_{\alpha}:\matK \longrightarrow
\matR^{n} \\
\hspace{0.8cm} x\longmapsto
\left(\sqrt{\alpha_{1}}\sigma_1(x),\ldots,\sqrt{\alpha_{n}}\sigma_{n}(x)\right)\end{array}
\right.$$ is called {\it twisted homomorphism}. When $\alpha=1$ the
twisted homomorphism is the {\it Minkowski homomorphism}.

It can be shown that if $I\subseteq {\matK}$ is a free
$\matZ$-module of rank $n$ with $\matZ$-basis $\{w_1,\ldots,w_n\}$,
then the image $\Lambda = \sigma_{\alpha}(I)$  is a lattice in
$\matR^n$ with basis $\{{ \sigma_{\alpha}(w_1)},\ldots,{
\sigma_{\alpha}(w_n)}\},$ or equivalently with generator matrix
${\bm M}= (\sigma_{\alpha}(w_{ij}))_{i,j=1}^{n}$ where $w_i =
(w_{i1,\cdots,w_{in}})$ for all $i=1,\cdots,n$.  Moreover, if
$\alpha I \overline{I} \subseteq \Delta(\matK|\matQ)^{-1}$ where
$\overline{I}$ denote the complex conjugation of $I,$ then
$\sigma_{\alpha}(I)$ is an integer lattice. Since $\matK$ is totally
real, the associated Gram matrix of $\sigma_{\alpha}(I)$ is ${\bm
G}={\bm M}.{\bm M^{t}} = \left( Tr_{\mattK|\mattQ}(\alpha
w_{i}\overline{w_{j}}) \right)_{i,j=1}^{n}$ \cite{Oggier}.

\begin{proposition}\label{detb}\cite{eva2} If $I\subseteq \matK$ is a fractional ideal, then for $\Lambda=\sigma_{\alpha}(I)$ and $det(\Lambda)=det(G)$, we have:
\begin{equation} \label{detb1} det(\Lambda)= det(G) =  N(I)^{2}
N_{\mattK|\mattQ}(\alpha)|d_{\mattK}|.\end{equation}
\end{proposition}

\begin{proposition} \label{distanciaminima}\cite{Oggier} Let $\matK$ be a totally real field number with  $[\matK:\matQ]=n$ and $I \subseteq
{\matK}$  a fractional ideal. The {\it minimum  product distance} of
$\Lambda =\sigma_{\alpha}(I)$ is
\begin{equation} d_{p,min}(\Lambda)= \sqrt{N_{\mattK|\mattQ}(\alpha)}min_{0\neq
y\in I}|N_{\mattK|\mattQ}(y)|.\end{equation} In particular, if $I$
is a principal ideal then $ d_{p,min}(\Lambda)=
\sqrt{\frac{det(\Lambda)}{|d_{\mattK}|}}.$
\end{proposition}

\begin{definition} \label{relative} The {\it relative minimum product distance} of $\Lambda,$ denoted by ${\bm
d_{p,rel}(\Lambda)}$, is the minimum product distance of a scaled
version of $\Lambda$ with unitary minimum norm vector.
\end{definition}

\section{Rotated $D_n$-lattices for $n=2^{r-2}$, $r \geq 5$ via
$\matK=\matQ(\zeta_{2^{r}}+\zeta_{2^{r}}^{-1})$}

In this section we will present some families of rotated
$D_n$-lattices using ideals and modules in the totally real number
field $\matK=\matQ(\zeta_{2^{r}}+\zeta_{2^{r}}^{-1}).$ One of the
strategies to construct these lattices was to start from the
standard characterization of $D_n$ as generated by the basis
\begin{equation}\label{beta} \beta= \{(-1,-1,0,\cdots,0),(1, -1,0,
\cdots,0),\cdots,(0,0,\cdots,1,-1)\}. \end{equation} We derive in
4.2 a rotated $D_n$-lattice as a sublattice of the rotated
$\matZ^{n}$ algebraic constructions presented in \cite{and},
\cite{Oggier} and \cite{gabriele}. Another strategy explored next in
4.1 is to investigate the necessary condition given in Proposition
\ref{detb}, for the existence of rotated $D_n$-lattices.

Let $\zeta=\zeta_{2^r}$ be a primitive $2^r$-th root of unity, $m =
2^{r},$ $\matK=\matQ (\zeta +\zeta^{-1})$ and
$n=[\matK:\matQ]=2^{r-2}$.

\subsection{\bf A first construction:} Let $\alpha \in \mathcal O_{\mattK}$ and $I \subseteq \mathcal O_{\mattK}$ an ideal. If $\sigma_{\alpha}(I)$ is a
rotated $D_n$-lattice scaled by $\sqrt{c}$, then
$det(\sigma_{\alpha}(I))=4\, c^n$. Based on Proposition \ref{detb},
taking $I=\mathcal O_{\mattK}$ and $c=2^{r-1}$, since
$d_{\mattK}=2^{(r-1)2^{r-2}-1}$ and $n=2^{r-2}$ it follows that a
necessary condition to construct a rotated $D_n$-lattice
$\sigma_{\alpha}(I)$ is to find an element $\alpha \in \mathcal
O_{\mattK}$ such that $N(\alpha)=8$. Table $1$ shows some elements
$\alpha \in \mathcal O_{\mattK}$ such that $N(\alpha)=8$ in low
dimensions. From it we got the suggestion for a general expression
for $\alpha$  as \begin{equation} \label{alpha} \alpha=4 +
(\zeta_{2^r}+ {\zeta{_{2^r}^{-1}}})- 2({\zeta{_{2^r}^{2}}}+
{\zeta{_{2^r}^{-2}}})-({\zeta{_{2^r}^{3}}}+ {\zeta{_{2^r}^{-3}}})
\end{equation} and then derive Proposition \ref{Prop1}.

\begin{table}[h]
\begin{center}
\begin{tabular}{c|c|c} \hline \hline
$r $& $\alpha $& ${N}(\alpha) $ \\ \hline  $4$ & $4 + (\zeta_{16}+
{\zeta{_{16}^{-1}}})- 2({\zeta{_{16}^{2}}}+
{\zeta{_{16}^{-2}}})-({\zeta{_{16}^{3}}}+ {\zeta{_{16}^{-3}}})$& $8$
\\ $5 $& $4 + (\zeta_{32}+ {\zeta{_{32}^{-1}}})-
2({\zeta{_{32}^{2}}}+ {\zeta{_{32}^{-2}}})-({\zeta{_{32}^{3}}}+
{\zeta{_{32}^{-3}}})$& $8$ \\  $6 $& $4 + (\zeta_{64}+
{\zeta{_{64}^{-1}}})- 2({\zeta{_{64}^{2}}}+
{\zeta{_{64}^{-2}}})-({\zeta{_{64}^{3}}}+ {\zeta{_{64}^{-3}}})$& $8$
\\  \hline \hline \end{tabular} \vspace{0.3cm} \caption{ }
\end{center}
\end{table}

To prove that $\frac{1}{\sqrt{2^{r-1}}}\sigma_{\alpha}(I)$ is a
rotated $D_n$-lattice we need the next  preliminary results.

\begin{proposition}\label{traco} \cite{and} If $\zeta=\zeta_{2^r}$  and $\matK=\matQ (\zeta +\zeta^{-1})$, then \newline $
Tr_{\mattK|\mattQ}(\zeta^k + \zeta^{-k})=\left\{\begin{array}{l} 0,
\ \
se \  \ gcd(k,2^r)<2^{r-1}; \\ -2^{r-1}, \ \ se \  \ gcd(k,2^r)=2^{r-1}; \\
2^{r-1}, \ \ se \  \ gcd(k,2^r)=2^{r}. \end{array}\right
.$\end{proposition}

\begin{proposition}\label{quase} If  $\matK=\matQ(\zeta+\zeta^{-1})$, $e_0 =1$ and $e_i=\zeta^{i}+\zeta^{-i}$ for $i = 1,\cdots,
{2^{r-2}-1}$, then \\
(a) $Tr_{\mattK|\mattQ}(\alpha e_i e_i) = \left\{\begin{array}{ll} 2^{r}, & \mbox{if }\, i=0,1 \\
2^{r+1}, & \mbox{if }\, 2 \leq i < 2^{r-2}-1  \\ 3.2^{r}, & \mbox{if
}\, i=2^{r-2}-1
\end{array} \right.$ \\
(b) $Tr_{\mattK|\mattQ}(\alpha e_i e_0) = \left\{\begin{array}{ll} 2^{r-1}, & \mbox{if }\, i=1 \\
-2^{r}, & \mbox{if }\, i=2 \\ -2^{r-1}, & \mbox{ if }\,\, i=3
\\  0, & \mbox{if }\, 3<i \leq 2^{r-2}-1 \end{array}\right.$ \\ (c) If $0 < i < j \leq 2^{r-2}-1$ then $$Tr_{\mattK|\mattQ}(\alpha e_i e_j) =
\left\{\begin{array}{ll} 2^{r-1}, & \mbox{ if }\,\, |i-j|=1\,\, \mbox{ and} \\
 & (i,j)\not\in \{(1,2),(2,1),(2^{r-2}-2, 2^{r-2}-1), \\ &
(2^{r-2}-1,2^{r-2}-2)\} \\
-2^{r}, & \mbox{ if }\,\, |i-j|=2 \\ -2^{r-1}, & \mbox{ if }\,\,
|i-j|=3  \\ 2^{r}, & \mbox{ if }\,\, (i,j) \in \{(2^{r-2}-2,
2^{r-2}-1), \\ & (2^{r-2}-1,2^{r-2}-2)\}
\\  0, & \mbox{ otherwise.}
\end{array}\right.$$
\end{proposition}
{\it Proof:} The proof is straightforward by calculating the
$gcd(k,2^{r})$ for some values of $k$ and applying Proposition
\ref{traco}. For $0 < i < j \leq 2^{r-2}-1$ we have:
\[\begin{split} Tr(\alpha e_i e_i) & = Tr(8)+ 4Tr({\zeta^{2i}}+
{\zeta^{-2i}})+ 2Tr({\zeta}+ {\zeta^{-1}})+
 Tr({\zeta^{2i+1}} + {\zeta^{-2i-1}}) \\& + Tr({\zeta^{2i-1}}+
{\zeta^{-(2i-1)}})- 4Tr({\zeta^{2}}+ {\zeta^{-2}})
-2Tr({\zeta^{2i+2}}+ {\zeta^{-(2i+2)}}) \\& -2Tr({\zeta^{2i-2}}+
{\zeta^{-(2i-2)}})-2Tr({\zeta^{3}}+ {\zeta^{-3}}) \\& -
Tr({\zeta^{2i+3}}+ {\zeta^{-(2i+3)}})  - Tr({\zeta^{2i-3}}+
{\zeta^{-(2i-3)}}) \end{split} \]  For $2\leq i < 2^{r-2}-1$ since
$gdc(k,2^{r})< 2^{r-1}$ for $k= 2i,2i \pm 1,2i\pm2, 2i \pm 3$ it
follows that $Tr(\alpha e_i e_i) = 2^{r+1}$. For $i=1,2^{r-2}-1$ the
development is analogous. For $i=0$ we have: \[ Tr(\alpha
e_0e_0)=Tr(4)+ Tr(\zeta+ {\zeta^{-1}}) -2Tr({\zeta^{2}}+
{\zeta^{-2}})-Tr({\zeta^{3}}+{\zeta^{-3}})=2^r\] and then it follows
(a).

 \[ \begin{split}  Tr(\alpha e_ie_0) &
=4Tr({\zeta^{i}}+{\zeta^{-i}})+Tr({\zeta^{i+1}}+{\zeta^{-(i+1)}}) +
Tr({\zeta^{i-1}}+{\zeta^{-(i-1)}}) \\&-2Tr({\zeta^{i+2}}+
{\zeta^{-(i+2)}})-2Tr({\zeta^{i-2}}+{\zeta^{-(i-2)}})\\&-
Tr({\zeta^{i+3}}+{\zeta^{-(i+3)}})-Tr({\zeta^{i-3}}+
{\zeta^{-(i-3)}})\end{split} \] For $i\neq 1,2,3$, since
$gcd(k,2^{r})< 2^{r-1}$ for $k=i, i\pm 1, i\pm 2, i\pm 3$ then
$Tr(\alpha e_ie_0)=0$.\\ For $i=1,2,3$ using $Tr(\zeta^{0}+
\zeta^{0})=2^{r-1}$ it follows (b).

 \[ \begin{split} Tr(\alpha e_ie_j)&
=Tr(\zeta^{i-j+1}+\zeta^{-(i-j+1)})+
Tr(\zeta^{i-j-1}+\zeta^{-(i-j-1)})
\\&-2[Tr(\zeta^{i-j+2}+\zeta^{-(i-j+2)})+
Tr(\zeta^{i-j-2}+\zeta^{-(i-j-2)})]\\& -[Tr(\zeta^{i-j+3}+
\zeta^{-(i-j+3)}) + Tr(\zeta^{i-j-3}+\zeta^{-(i-j-3)})]\\& -2
Tr(\zeta^{i+j-2}+\zeta^{-(i+j-2)})-
Tr(\zeta^{i+j+3}+\zeta^{-(i+j+3)})\\&-
Tr(\zeta^{i+j-3}+\zeta^{-(i+j-3)}) \end{split}\]

Since $gcd(k,2^{r})<2^{r-1}$ for $k=i+j+1, i+j+ 2$;
$gcd(i+j+3,2^{r})<2^{r-1}$ for $i+j \neq 2^{r-1}-3$; $gcd(i+j+3,2^r)
=2^{r-1}$, for $i+j = 2^{r-1}-3$ and $gcd(i+j-3,2^r)<2^{r-1}$, for
$i+j \neq 3$, it follows (c).

\bb

\begin{proposition}\label{Prop1} The lattice  $\frac{1}{\sqrt{2^{r-1}}}\sigma_{\alpha}(\mathcal
O_{\mattK}) \subseteq \matR^{2^{r-2}}$, $\alpha =4 + (\zeta_{2^r}+
{\zeta{_{2^r}^{-1}}})- 2({\zeta{_{2^r}^{2}}}+
{\zeta{_{2^r}^{-2}}})-({\zeta{_{2^r}^{3}}}+ {\zeta{_{2^r}^{-3}}})$
is a rotated $D_n$-lattice for $n=2^{r-2}$.
\end{proposition}
{\it Proof:} The Gram matrix for
$\frac{1}{\sqrt{2^{r-1}}}\sigma_{\alpha}(\mathcal O_{\mattK})$
related to the $\matZ$-basis  $\{e_0,e_1,\cdots,e_{n-1}\}$ is
\[{\bm G}={ \left(
\begin{array}{ccccccccccc}
 2 & 1 & -2 & -1 & 0 & \cdots  &  &  &  & \cdots  & 0 \\
 1 & 2 & 0 & -2 & -1 & 0 &  &  &  &  & \vdots  \\
\vdots  & 0 & \ddots  & \ddots & \ddots & \ddots  & \ddots  & \ddots  & \ddots  & 0 & \vdots  \\
  &  &  &  & \ddots  & 0 & -1 & -2 & 1 & 4 & 2 \\
 0 & \vdots &  &  &  & \cdots & 0 & -1 & -2 & 2 & 6
\end{array}
\right)}\] \noindent and it is easy to see that ${\bm G}$ is the
Gram matrix for $D_n$ related to the generator matrix ${\bm T}{\bm
B}$ where
\[{\bm T}={\left(
\begin{array}{cccccccc}
 0 & \cdots &  &  &  & \cdots & 0 & -1 \\
 0 & \cdots &  &  & \cdots  & 0 & 1 & 0 \\
 0 & \cdots &  & \cdots & 0 & -1 & 0 & 1 \\
 \vdots & \vdots  & \vdots & \vdots  & \vdots  &  \vdots & \vdots  &  \\
 0 & 0 & 1 & 0 & -1 & 0 & \cdots & 0 \\
 0 & -1 & 0 & 1 & 0 & \cdots & \cdots  & 0 \\
 1 & -1 & -1 & 0 & \cdots &  & \cdots  & 0
\end{array}
\right)}\] and ${\bm B}$ is the standard generator matrix for $D_n$
given by basis $\beta$ (\ref{beta}). So, since lattices with the
same Gram matrix must be Euclidean equivalent, then  then
$\sigma_{\alpha}(I)$ is a rotated $D_n$-lattice. \bb

We determine next the relative minimum product distance of the
rotated $D_n$-lattice considered in Proposition \ref{Prop1}.

Using Propositions \ref{detb} and \ref{Prop1} we conclude:

\begin{corollary} If $m = 2^{r}$, $r \geq 4$, $\matK=\matQ(\zeta_m+\zeta_m^{-1})$ and $\alpha =4 + (\zeta_{2^r}+ {\zeta{_{2^r}^{-1}}})-
2({\zeta{_{2^r}^{2}}}+ {\zeta{_{2^r}^{-2}}})-({\zeta{_{2^r}^{3}}}+
{\zeta{_{2^r}^{-3}}})$ then $N_{\mattK|\mattQ}(\alpha) = 8.$
\end{corollary}

\begin{proposition} For $n=2^{r-2}$, if $\Lambda = \frac{1}{\sqrt{2^{r-1}}}\sigma_{\alpha}(\mathcal
O_{\mattK})$ and $\alpha$ as in (\ref{alpha}) then the lattice
relative minimum product distance is

$${\bm
d_{p,rel}\left(\frac{1}{\sqrt{2^{r-1}}}\sigma_{\alpha}(\mathcal
O_{\mattK})\right)} =2^{\frac{3-r n}{2}}.$$
\end{proposition} {\it Proof:} The minimum norm of the standard $D_n$ is
$\sqrt{2}.$  $\mathcal O_{\mattK}$ is a  principal ideal, therefore
using Proposition \ref{distanciaminima} we have
$d_{p,min}(\sigma_{\alpha}(\mathcal O_{\mattK})) =
\sqrt{N(\alpha)N(\mathcal O_{\mattK})^{2}}.$ Since $N(\alpha)=8$ and
$N(\mathcal O_{\mattK})=1,$ then

$${\bm
d_{p,rel}\left(\frac{1}{\sqrt{2^{r-1}}}\sigma_{\alpha}(\mathcal
O_{\mattK})\right)} =
\frac{1}{\sqrt{2}^{n}}\frac{1}{\sqrt{2^{r-1}}^{n}} \sqrt{8}
=\frac{\sqrt{8}}{2^{r \frac{n}{2}}} = 2^{\frac{3-rn}{2}}.$$ \bb
\vspace{1cm}

\subsection{\bf A second construction:} In \cite{and} and \cite{gabriele}
families of rotated $\matZ^n$-lattices obtained as image of a
twisted homomorphism applied to $\matZ[\zeta+\zeta^{-1}]$ and having
full diversity are constructed. Those constructions consider $\alpha
= 2 + e_1$ and $\alpha=2-e_1,$ respectively, and generate equivalent
lattices in the Euclidean metric by permutations and coordinate
signal changes.

We will use in our construction the rotated $\matZ^{n}$-lattice
 $\Lambda=\frac{1}{\sqrt{2^{r-1}}}\sigma_{\alpha}(I)$ with
$\alpha = 2 + e_1$ and $I=\mathcal
O_{\mattK}=\matZ[\zeta+\zeta^{-1}],$ and then consider $D_n$ as a
sublattice of $\Lambda$.

If $e_0 =1$ and $e_i=\zeta^{i}+\zeta^{-i}$ for $i = 1,\cdots,
{2^{r-2}-1}$, by \cite{gabriele} a generator matrix for the rotated
$\matZ^{n}$-lattice $\Lambda=$
$\frac{1}{\sqrt{2^{r-1}}}\sigma_{\alpha}(\mathcal O_{\mattK})$ is
${\bm M_1} = \displaystyle\frac{1}{\sqrt{2^{r-1}}}{\bm N} {\bm A},$
where ${\bm N}=(\sigma_i(e_{j-1})_{i,j=1}^{n}$ and ${\bm A}=
diag(\sqrt{\sigma_{k}(\alpha)}).$ Let ${\bm T}$ the  basis change
matrix  \[{{\bm  T}=\left( \begin{array}{ccccc} 1 & -1  & \cdots & 1 & -1 \\
1 & -1 & \cdots  & 1 & 0 \\ \vdots & \vdots & \ddots  & \vdots & \vdots \\
1  & 0   & \cdots  & 0 & 0
\end{array}\right)}.\]
For ${\bm M} = {\bm T}{\bm M_1}$, ${\bm G}={\bm M} {\bm M^{t}}={\bm
I_n}$ and we will consider the standard lattice $D_n \subseteq
\matZ^{n}$ rotated by ${\bm M}$.

\begin{proposition}\label{idealI} Let $I \subseteq \mathcal O_{\mattK}$ be the $\matZ$-module with
$\matZ$-basis \linebreak $\{-e_1,e_2,\cdots,-e_{n-1},  -2e_0 + 2 e_1
- 2e_2 + \cdots -2e_{n-2} + e_{n-1}\}$ and $\alpha=2+e_1.$ The
lattice $\frac{1}{\sqrt{2^{r-1}}}\sigma_{\alpha}(I) \subseteq
\matR^{2^{r-2}}$ is a rotated $D_n$-lattice.
\end{proposition} \noindent{\it Proof:} Let ${\bm B}$ be the generator matrix of $D_n$ associated to the basis
$\beta$ (\ref{beta}).  Using homomorphism properties, a
straightforward computation shows that ${\bm B} {\bm M} =$

$$ {
 \frac{1}{\sqrt{2^{r-1}}}\footnotesize\left(\begin{array}{ccc} \sigma_1(-2 e_0 +
 \cdots - 2e_{n-2} + e_{n-1} ) & \cdots & \sigma_n(-2 e_0 +
 \cdots - 2e_{n-2} + e_{n-1} )
\\ \sigma_1(-e_{n-1}) & \cdots &
\sigma_n( -e_{n-1})
\\ \vdots & \ddots & \vdots \\ \sigma_1(-e_1)
& \cdots & \sigma_n(-e_1) \end{array}\right) {\bm A}} $$ is a
generator matrix for $\frac{1}{\sqrt{2^{r-1}}}\sigma_{\alpha}(I).$
This lattice is a rotated $D_n$-lattice since ${\bm B} {\bm M} ({\bm
B} {\bm M})^{t} = {\bm B} {\bm B^{t}}$ is the standard Gram matrix
of $D_n$ relative to the basis $\beta$. \bb

We show next that the rotated $D_n$-lattice of the last proposition
is associated to a principal ideal of $\mathcal O_{\mattK}$ and then
calculate its relative  minimum product distance.

\begin{proposition} Let $I$ be the $\matZ$-module given in  Proposition \ref{idealI}.
$I$ is a principal ideal  and $I =  e_1 \mathcal O_{\mattK}.$
\end{proposition} {\bf Proof:} It is easy to see that $I =  2e_0\matZ +  e_1\matZ + \cdots +
e_{n-1}\matZ.$ Let $x \in e_1 \mathcal O_{\mattK}.$ Then $x = e_1
(a_0 e_0 + a_1 e_1 + a_2 e_2 + \cdots + a_{n-1}e_{n-1}) =   a_0 (e_1
+ e_{-1}) + a_1 (e_2 + 2e_0) + a_2 (e_3 + e_{-1}) + \cdots +
a_{n-1}(e_n + e_{-n+2}) = a_1 (2 e_0) + (2 a_0 + a_2)(e_1) + (a_1
+a_3)(e_2) + \cdots + (a_{n-2})(e_{n-1}) \in I.$ Now, if $x \in I,$
then  $x = a_0 2e_0 + a_1 e_1 + \cdots + a_{n-1} e_{n-1} = (e_1)[a_0
e_1 + a_1 e_2 + (a_2 - a_0) e_3  + ( a_3 - a_1) e_4 + (a_4 - a_2
-a_0) e_5 + (a_5 - a_3 -a_1)e_6 + \cdots + ( a_{n-1})e_{n-2} +
(a_{n-2} - a_{n-4} \cdots - a_0)e_{n-1}] \in e_1 \mathcal
O_{\mattK}.$ So, $I$ is a principal ideal of $\mathcal O_{\mattK}.$
\bb

\begin{remark} It follows from Proposition \ref{distanciaminima} and
Definition \ref{relative} that the \ relative minimum product
distance of $D_n$-lattices constructed from principal ideals in
$\mathcal O_{\mattK}=\matQ(\zeta_m+\zeta_m^{-1})$, $m=2^{r}$, $r
\geq5$, depends only of the determinant of $D_n$ and of the
discriminant of $\matK$. Therefore for any  construction of a
rotated $D_n$ lattice from a principal ideal $I$ in $\mathcal
O_{\mattK}$ the relative minimum product distance is
$d_{p,rel}(\sigma_{\alpha}(I)) = 2^{\frac{3-r n}{2}}.$
\end{remark}

It is also interesting to note that besides being Euclidean
equivalent, the lattices obtained through the first and second
constructions are equivalent in the Lee metric since the isometry is
a composition of permutations and coordinate signal changes.

The density $\Delta(\Lambda)$ of a lattice $\Lambda \subseteq
\matR^{n}$ is given by $\Delta(\Lambda) =
\frac{(d/2)^{n}Vol(B(1))}{det(\Lambda)^{1/2}}$ where $Vol(B(1))$ is
the volume of the unitary sphere in $\matR^{n}$ and $d$ is the
minimum norm of $\Lambda$. The parameter $\delta(\Lambda) =
\frac{(d/2)^{n}}{det(\Lambda)^{1/2}}$ is so called center density.
Table $2$ shows a comparison between the normalized $d_{p,rel}$ and
the center density of rotated $\matZ^n$-lattices constructed in
\cite{gabriele} and rotated $D_n$-lattices constructed here via
principal ideals in $\matK=\matQ(\zeta+\zeta^{-1})$, $n=2^{r-2}$.
Asymptotically we have
\begin{equation} \label{limite} \lim_{n\longrightarrow
\infty}\frac{\sqrt[n]{d_{p,rel}(\matZ^{n})}}{\sqrt[n]{d_{p,rel}(D_n)}}
=\sqrt{2}\,\, \, \mbox{and}\,\,\,  \lim_{n\longrightarrow
\infty}\frac{\delta(\matZ^{n})}{\delta(D_n)} =0. \end{equation}

\begin{table}[h]
\begin{center}
\begin{tabular}{|c|c|c|c|c|c|}
\hline \hline
$r$ & $n$ & $\sqrt[n]{d_{p,rel}(\matZ^n)}$ & $\sqrt[n]{d_{p,rel}(D_n)}$ & $\delta(\matZ^n)$ & $\delta(D_n)$ \\
\hline
$4$ & $4$ & $ 0.385553  $ & $ 0.324210 $ &$0.062500$ & $0.125000$\\
\hline
$5$ & $8$ & $ 0.261068 $ & $ 0.201311 $ &$0.003906$ & $0.031250$ \\
\hline
$6$ & $16$ & $ 0.180648 $ & $ 0.133393 $ &$0.000015$ & $0.001953$ \\
\hline
$7$ & $32$ & $ 0.126361 $ & $ 0.091307 $ &$2.3 \times 10^{-10}$ & $7.6 \times  10^{-6}$\\
\hline
$8$ & $64$ & $ 0.088868 $ & $ 0.063523  $ &$5.4 \times  10^{-20}$ & $1.1 \times  10^{-10}$\\
\hline
$9$ & $128$ & $ 0.062669 $ & $ 0.044554 $ &$2.9 \times  10^{-39}$ & $2.7 \times  10^{-20}$\\
\hline \hline
\end{tabular}
\vspace{0.3cm} \caption{ }
\end{center}
\end{table}

If the  goal is to construct lattices which have good performance on
both Gaussian and Rayleigh channels, we may assert that taking into
account the trade-off density versus product distance, there is some
advantages in considering these rotated $D_n$-lattices instead of
rotated $\matZ^{n}$-lattices, $n=2^{r-2}$, $r \geq 5,$ in high
dimensions.

\section{Rotated $D_n$-lattices for $n=\frac{p-1}{2}$, $p$ prime,
via $\matK=\matQ(\zeta_{p}+\zeta_{p}^{-1})$}

Let $\zeta=\zeta_{p}$ be a primitive $p$-th root of unity, $p$
prime,  $\matL=\matQ (\zeta)$ and $\matK=\matQ (\zeta+\zeta^{-1})$.
We will construct a family of rotated $D_n$-lattices, derived from
the construction of a rotated $\matZ^{n}$-lattice in \cite{Oggier},
via a $\matZ$-module that is not an ideal. Let
$e_j=\zeta^{j}+\zeta^{-j}$ for $j = 1,\cdots, {(p-1)/2}.$

By \cite{Oggier} a generator matrix of the rotated
$\matZ^{n}$-lattice
$\Lambda=\frac{1}{\sqrt{p}}\sigma_{\alpha}(\mathcal O_{\mattK})$ is
${\bm M}= \displaystyle\frac{1}{\sqrt{p}}{\bm T} {\bm N} {\bm A} ,
\mbox{ where }$ ${\bm T}=(t_{ij})$ is an upper triangular matrix
with $t_{ij}=1$ if $i\leq j,$ ${\bm N} =
(\sigma_i(e_j))_{i,j=1}^{n}$ and ${\bm A}
=diag(\sqrt{\sigma_{k}(\alpha)}).$ We have ${\bm G}={\bm M} {\bm
M^{t}}={\bm I_n}$ \cite{Oggier}.

\begin{proposition} \label{rotacionado} Let $I \subseteq \mathcal O_{\mattK}$ be a $\matZ$-module with
$\matZ$-basis  $\{e_1,e_2,\cdots,e_{n-1},-e_1 - 2e_2 - \cdots
-2e_n\}$ and $\alpha=2-e_1.$ The lattice
$\frac{1}{\sqrt{p}}\sigma_{\alpha}(I) \subseteq
\matR^{\frac{p-1}{2}}$ is a rotated $D_n$-lattice.
\end{proposition} \noindent{\it Proof:} Let ${\bm B}$ be a generator matrix for $D_n$ given by basis $\beta$ \ref{beta}.
Using homomorphism properties, a straightforward computation shows
that ${\bm B}{\bm M}$ is a generator matrix for
$\Lambda=\frac{1}{\sqrt{p}}\sigma_{\alpha}(I).$ This lattice is a
rotated $D_n$ since ${\bm B} {\bm M} ({\bm B} {\bm M})^{t} = {\bm B}
{\bm B^{t}}$ is a Gram matrix of $D_n.$ It has full diversity since
it is contained in $\frac{1}{\sqrt{p}}\sigma_{\alpha}(\mathcal
O_{\mattK})$ \cite{Oggier}. \bb

\begin{proposition} The $\matZ$-module $I \subseteq \mathcal O_{\mattK}$ is not an ideal of $\mathcal O_{\mattK}$. \end{proposition}
\noindent{\it Proof:} The set $\{e_1,e_2,\cdots,e_{n-1},2e_n\}$ is
an another $\matZ$-basis to $I$. We will show that $e_n$ is not in
$I.$ Indeed, if $e_n \in I$, then $I = \mathcal O_{\mattK},$ but
$\left|\frac{\mathcal O_{\mattK}}{I}\right| = 2.$ So, $e_n \not\in
I.$  $e_{n-1} e_{1}$ is not in $I$. In fact, note that $e_{n-1} e_1
= e_{n} + e_{n-2}$ and $e_{n-2} \in I$. If $e_{n-1}e_1 \in I$, then
$e_n = e_{n-1}e_1 - e_{n-2} \in I,$ and this doesn't happen. \bb

\begin{proposition} If $\Lambda =\frac{1}{\sqrt{p}}\sigma_{\alpha}(I) \subseteq
\matR^{\frac{p-1}{2}}$ with $\alpha$ and $I$ as in the Proposition
\ref{rotacionado}, then the relative minimum product distance is

$${\bm d_{p,rel}}(\Lambda)= 2^{\frac{1-p}{4}} p^{\frac{3-p}{4}} .$$
\end{proposition} \noindent{\it Proof:} First  note that $|N(e_1)| =1.$
Indeed, $(\zeta+\zeta^{-1})(-\zeta^{p-1} - \zeta^{p-2} - \cdots -
\zeta - 1)=1$ and so

$$N(\zeta+\zeta^{-1})N(-\zeta^{p-1} -
\zeta^{p-2} - \cdots - \zeta - 1)=N(1)=1.$$

Since $e_1 \in \mathcal O_{\mattK}$, then $N(e_1) \in \matZ$, what
implies $|N(e_1)|=1.$ Now, the minimum norm in $D_n$ is $\sqrt{2}.$
By Proposition \ref{distanciaminima},  ${\bm
d_{p}(\sigma_{\alpha}(I))} =\sqrt{N(\alpha)} min_{0\neq y \in
I}|N(y)| = \sqrt{p},$ since $min_{0\neq y \in I}|N(y)| =1.$
Therefore, the relative minimum product distance is

$$d_{p,rel}\left(\frac{1}{\sqrt{p}}\sigma_{\alpha}(I)\right) =
\left(\frac{1}{\sqrt{p}^{\frac{p-1}{2}}}\right)
\left(\frac{1}{\sqrt{2}^{\frac{p-1}{2}}}\right)\sqrt{p} =
2^{\frac{1-p}{4}} p^{\frac{3-p}{4}}.$$ \bb \vspace{0.4cm}

Table $3$ shows a comparison between the normalized $d_{p,rel}$ and
the center density $\delta$ of rotated $\matZ^n$-lattices
constructed in \cite{Oggier} and rotated $D_n$-lattices constructed
here, $n=(p-1)/2$. As in Section \ref{limite} we also have for
$\Lambda = \frac{1}{\sqrt{p}}(\sigma_{\alpha}(I)) \subseteq
\matR^{\frac{p-1}{2}}$ and $p$ prime, the following results:
$$\lim_{n\longrightarrow
\infty}\frac{\sqrt[n]{d_{p,rel}(\matZ^{n})}}{\sqrt[n]{d_{p,rel}(D_n)}}
=\sqrt{2}\,\,  \mbox{and} \,\, \lim_{n\longrightarrow
\infty}\frac{\delta(\matZ^{n})}{\delta(D_n)} =0.$$

\begin{table}[h]
\begin{center}
\begin{tabular}{|c|c|c|c|c|c|}
\hline \hline
$p$ & $n$ & $\sqrt[n]{d_{p,rel}(\matZ^n)}$ & $\sqrt[n]{d_{p,rel}(D_n)}$ & $\delta(\matZ^n)$ & $\delta(D_n)$ \\
\hline
$11$ & $5$ & $ 0,38321 $ & $0,27097 $ & $0,03125$ & $0,08838$ \\
\hline
$13$ & $6$ & $ 0,34344 $ & $0,24285 $ & $0,01563$ & $0,06250$ \\
\hline
$17$ & $8$ & $ 0,28952$ & $0,20472 $  & $0,00390$ & $0,03125$ \\
\hline
$19$ & $9$ & $0,27187 $ & $ 0,19105 $  & $0,00195$ & $0,02209$ \\
\hline
$23$ & $11$ & $ 0,24045$ & $ 0,17003 $  & $0,00049 $ & $0,01105$\\
\hline \hline
\end{tabular}
\vspace{0.3cm} \caption{ }
\end{center}
\end{table}

\section{Conclusion}

In this work we construct some families of full diversity rotated
$D_n$-lattices. These lattices present good performance for signal
transmission over both Gaussian and Rayleigh channels. Considering
the trade-off between density and relative product distance we may
assert that the rotated $D_n$-lattices presented here have better
performance than the known rotated $\matZ^{n}$-lattices for
$n=2^{r-2},$ $r \geq 5$ and $n=\frac{p-1}{2}$, $p$ prime, $p \geq
7$.

\end{document}